# Hybrid Numerical Modeling of Ballistic Clay under Low-Speed Impact using Artificial Neural Networks

YeonSu Kim[a, b], Yoon A Kim[a], Seo Hwee Park[a], YunHo Kim[a]*

[a]*Extreme Environments and Impact Lab., School of Mechanical Engineering, Pusan National University, Busan, Republic of Korea*

[b]*Department of Aerospace Engineering, Pusan National University, Busan, Republic of Korea*

*Corresponding author

Email address

YeonSu Kim : lelghan@pusan.ac.kr

Yoon A Kim : rladbsdk0825@pusan.ac.kr

Seo Hwee Park : suhwee12@pusan.ac.kr

YunHo Kim : yunho@pusan.ac.kr

**Abstract**

Roma Plastilina No. 1 clay has been widely used as a conservative boundary condition in bulletproof vests, namely to play the role of a human body. Interestingly, the effect of this boundary condition on the ballistic performance of the vests is indiscernible. Moreover, back face deformation should be characterized by measuring the indentation in the deformed clay, which is important for determining the lethality of gunshots. Therefore, several studies have focused on modeling not only bulletproof vests but also the clay backing material. Despite various attempts to develop a suitable numerical model, determining the appropriate physical parameters that can capture the high-strain-rate behavior of clay is still challenging. In this study, we predicted indentation depth in clay using an artificial neural network (ANN) and determined the optimal material parameters required for a finite element method (FEM)-based model using an inverse tracking method. Our ANN–FEM hybrid model successfully optimized high-strain-rate material parameters without the need for any independent mechanical tests. The proposed novel model achieved a high prediction accuracy of over 98% referring impact cases.

# 1. Introduction

An artificial neural network (ANN) is a computational learning algorithm that is inspired by the biological neural network present in the human brain **(Dongare et al., 2012)**. ANNs are typically composed of multiple layers of neurons such that the input data pass repeatedly through these layers before the output data are generated **(Mishra and Srivastava, 2014)**. Owing to this distinctive learning mechanism, ANNs are expected to solve highly complex problems; however, such networks cannot distinguish between the cause and effect of the underlying problem. Therefore, "black-box" models such as ANNs are frequently arguable compared to conventional numerical approaches, such as the finite element method (FEM). Nevertheless, it is possible to incorporate ANNs into conventional methods. For example, various types of impacts, including large-scale deformations and complex failure modes, can be effectively simulated by combining these two approaches **(Fernández et al., 2008; Ryan et al., 2016; KılıÇ et al., 2015; Ramasamy and Sampathkumar, 2014; Mousavi and Khoramishad, 2019)**. For example, Fernández- Fernández et al. **(Fernández et al., 2018)** designed a multilayer perceptron model to predict the ballistic performance of carbon fiber-reinforced polymers under high-velocity oblique impacts. Ryan et al. **(Ryan et al., 2016)** used an ANN to predict the performance of multiwall aluminum Whipple shields against a wide range of hypervelocity impacts with a high accuracy of 92%. KılıÇ et al. **(KılıÇ et al., 2015)** employed a combined FEM–ANN method to predict the ballistic penetration depth of steel armors with high accuracy. Ramasamy and Sampathkumar **(Ramasamy and Sampathkumar, 2014)** evaluated the compressive strength of composites using ANNs.

ANNs require sufficient training data to make accurate predictions; otherwise, various numerical problems arise, such as overfitting. Overfitting implies that the prediction accuracy of the validation datasets is much lower than that of the training datasets **(Jordan and Mitchell,**

**2015)**. This can cause a significant loss in the prediction accuracy of the model. However, it is difficult to determine the optimum number of training datasets required by ANNs. In addition, it is challenging to acquire sufficient training data for problems such as predicting the impact of high-strain-rate deformation on different materials. Moreover, there are no set guidelines for determining the optimum number of experiments necessary, as well as the important physical outcomes for a given problem.

Roma Plastilina (RP) No. 1 clay is an oil-based backing material commonly used for evaluating the performance of body armors in ballistic tests as per the U.S. National Institute of Justice (NIJ) standards (**Council, 2009; Lehowicz, 2010; Council 2012; Standard, 1987; Standard, 2008**). In such impact tests, the extent of damage to the human body by a nonpenetrating gunshot wound is determined by measuring the depth of indentation in the RP clay. RP clay provides additional support to a target and is also expected to deteriorate significantly to ballistic performance. Therefore, it is important to investigate the effectiveness of RP clay as a backing material. Moreover, RP clay backing material needs to be validated before the start of every impact test in accordance with the NIJ standards (**Standard, 1987; Standard, 2008**). To this end, several studies have been conducted to model and characterize RP clay (**Buchely et al., 2016; Graham and Zhang, 2019; Hernandez et al., 2015**).

The FEM has been used to simulate RP clay (**Johnson, 1983**). Mates et al. (**Mates et al., 2014**) obtained the material parameters of the Johnson–Cook (J–C) model by using a reference strain rate of 0.118 $s^{-1}$ and reference temperature of 23 °C. Gad and Gao (**Gad and Gao, 2020**) compared indentation depth of between J-C model and new constitutive model, which can be applied in RP clay modeling based on a case where a 44.5mm cylindrical indenter drops from 2-m height referring to (**Standard, 2008**), to confirm effect of both temperature and strain rate. Hernandez et al. (**Hernandez et al., 2021**) used a model based on the J-C model and an inverse method to obtain the optimized set of material parameters characterizing RP clay including a test in which a 63.5mm spherical indenter drops from 2-m height. Gilson et al. (**Gilson et al., 2020**) analyzed the effect of Young's modulus (in the range 2–6 MPa) on the ballistic response of RP clay by comparing the results of numerical simulations and physical experiments correspond to where a 63.5mm spherical indenter dropped from 2-m heights. Nevertheless, the previous models have only been applied to a few limited cases involving 63.5 or 44.5mm diameter indenter from 2-m height. The model may not cover practical impact conditions of various geometry of indenter and wide range of strain rates. Despite of the importance of the clay, there are no practicable numerical models of the clay yet due to the difficulties of obtaining reliable material parameters in various loading conditions.

In this study, we used ANNs to determine the optimal values of the material parameters required by an FEM-based model for characterizing RP clay in accordance with the NIJ standards and their experimental results. First, we modeled the indenter and RP clay to obtain the necessary training datasets. Second, three ANNs were designed based on the results of the FEM simulations that could predict the indentation depths in RP clay due to impact by a spherical indenter with a diameter of (1) 44.5 mm at 4.47 m/s, (2) 44.5 mm at 6.26 m/s, and (3) 63.5 mm at 6.26 m/s. Next, we selected the optimal material parameters based on the predictions of our ANNs using an inverse tracking method. Finally, we implemented these optimal material parameters in our model for verification. We found that the accuracy of the FEM-based model was significantly improved and that the mean relative error was reduced from 17.27% to 1.74%.

## 2. Method

The procedure for determining the optimal material parameters used to model RP clay in this work is shown in **Figure. 1**. To validate the RP clay before use in bulletproof armors, several impact tests using an indenter were conducted following the NIJ standards. The diameters of the indenter were set to 44.5 mm and 63.5 mm in our numerical models **(Standard, 1987; Standard, 2008)**. In addition to these two cases, a case with a lower impact velocity was studied using the indenter of diameter 44.5 mm **(Kim, 2018)**. First, we numerically modeled these three impact cases (Case 1, 2, and 3). Next, we designed three independent ANN models ($A_1$, $B_2$, and $C_3$), which predicted indentation depths ($K_{1j}$, $K_{2j}$, and $K_{3j}$) corresponding to the three impact cases by learning from the FEM dataset as the material parameters ($A_j$) were varied. Subsequently, the predicted values of the material parameters were compared with the corresponding reference results represented by the function $f$. Finally, $f$ was minimized and the optimal material parameters (A) were determined using the inverse tracking method ($f^{-1}$).

### 2.1 Finite element modeling

The parameter settings and experimental results of the different impact cases are listed in **Table 1**. Case 1 corresponds to a cylindrical indenter made of 4340 steel with a spherical head of diameter 44.5±0.5 mm that was dropped from a height of 2 m **(Kim, 2018)**. Cases 2 and 3 correspond to two existing NIJ standards **(Standard, 1987; Standard, 2008)**. Case 1 was introduced to facilitate high-accuracy predictions by the ANN.

2D axisymmetric FEM simulations were carried out using the Ansys Autodyn software. Kim Y. A et al. **(Kim et al., 2022)** compared the indentation depth according to the effect of mesh size and boundary conditions of both clay and indenter for elaborate RP clay modeling. Appropriate mesh size and boundary conditions which are used in our numerical modeling are referred given in **Table. 2**. Moreover, both the indenter and RP clay were modeled as shell elements, as shown in **Figures. 2 and 3**. The duration of all the simulations was more than 10 ms, which was sufficiently long to slow down the indenter. At the end of each simulation, the displacement of the indenter was recorded and compared with referred indentation depth given in **Table 1**.

Modeling ballistic clay in Autodyn requires an EOS as well as the J–C model to account for both the hydrostatic (EOS) and deviatoric (J–C) components of the stress. For the indenter, the input parameters of 4340 steel were supplied from the Autodyn material library. For the RP clay, the parameters of a polynomial EOS and the J–C model were supplied based on previous results, which are listed in **Table 3 (Mates et al., 2014; Gad and Gao, 2020)**. Note that three different EOSs were considered; however, the choice of the EOS had a negligible effect on the RP clay. The EOS we adopted can be described by the following two equations:

$$p = A_1\mu + A_2\mu^2 + A_3\mu^3 + (B_0 + B_1\mu)p_0 e \qquad \mu > 0 \text{ (compression)} \qquad (1)$$

$$p = A_1\mu + A_2\mu^2 + B_0 p_0 e, \qquad \mu < 0 \text{ (tension)} \qquad (2)$$

where $A_1$, $A_2$, $A_3$, $B_0$, and $B_1$ represent the material parameters; $\mu$ is the compressibility; $p_0$ is the zero-pressure density; and $e$ is the internal energy per unit mass. The J–C model is given by the following equation:

$$\sigma = (A + B\varepsilon_p^n)\left(1 + C \ln\frac{\dot{\varepsilon}_p}{\dot{\varepsilon}_r}\right)\left[1 - \left(\frac{T-T_0}{T_m-T_0}\right)^m\right] \quad , \quad (3)$$

where *A*, *B*, *C*, *n*, and *m* are the material parameters, which represent the initial yield stress, hardening constant, strain rate constant, hardening exponent, and thermal softening exponent, respectively; $\varepsilon_p$ is the effective plastic strain; $\dot{\varepsilon}_p$ is the normalized effective plastic strain rate; $\dot{\varepsilon}_r$ is the reference strain rate; $T_0$ is the reference temperature; and $T_m$ is the melting temperature.

The instantaneous erosion strain (ISE) was applied to erode elements which have excessive level of strain. The material parameters in equations (1-3) were referred to the previous works **(Mates et al., 2014; Gad and Gao, 2020)**. A small value for the initial yield stress *A* was selected (i.e., 0.01 kPa) from previous studies to model the high plasticity of RP clay. Keeping the above considerations in mind, we considered only eight material parameters ($A_1$, $A_2$, *B*, *n*, *C*, $\varepsilon$, *m*, and ISE) in our final model that could be determined unambiguously. These are listed in **Table 4**.

## 2.2 Artificial neural network

An ANN typically consists of an input layer, a hidden layer, and an output layer, as shown in **Figure. 4**. The flow of data through the ANN can be described by the following equation (**McClelland et al., 1987**):

$$\underline{y_i} = \underline{w} \cdot x_i + b_i \quad , \tag{4}$$

where $x_i$ is the input data, $\underline{w}$ is the weight, and $b_i$ is the bias. The sum of $\underline{y_i}$ was transformed using an activation function (**Ramachandran et al., 2017; Sharma et al., 2017**). A rectified linear unit (ReLU) activation function was used between the input and hidden layers, as well as between the multiple hidden layers (**Nair and Hinton, 2010**). ReLU is defined as:

$$f(x) = max\ (0, x) \quad . \tag{5}$$

A linear activation function was used between the hidden and output layers, such that

$$f(x) = ax \quad . \tag{6}$$

An ANN learns iteratively by dividing the input data into training and validation data (**Smith, 2018**). In this study, we designed ANN models using Google Colab, which is a free Python environment. The training and validation data were split into an 80:20 ratio, and each model was trained the same number of times for 20,000 epochs. The adaptive moment estimation algorithm (ADAM), which is based on gradient descent, was used to optimize the performance of our model **(Kingma and Ba, 2015)**. The mean squared error (MSE) was used as the loss function, which is given by **(Allen, 1971).** In addition, we used the backpropagation algorithm to train our model **(Li et al., 2012)**.

$$MSE = \frac{1}{N}\sum_{i=0}^{n}(y_i - \hat{y}_i)^2 \ . \tag{7}$$

$(y_i: actual\ value, \hat{y}_i: predicted\ value)$

To find the ANN model with the best prediction accuracy, we first varied the total number of neurons (i.e., 25, 50, 75, and 100) based on one hidden layer. Next, we increased the number of hidden layers (i.e., 1, 2, and 3) to evaluate their effect on the performance of the ANN. Root mean square error (RMSE) and coefficient of determination ($R^2$) were used as the evaluation indices in this study, which are given by **(Wang and Lu, 2018; Nagelkerke, 1991)**.

$$RMSE = \sqrt{\frac{1}{N}\sum_{i=0}^{n}(y_i - \hat{y}_i)^2} \tag{8}$$

$(y_i: actual\ value, \hat{y}_i: predicted\ value)$

$$R^2 = \frac{\sum_{i=0}^{n}(\hat{y}_i - \underline{y})^2}{\sum_{i=0}^{n}(y_i - \underline{y})^2} = 1 - \frac{\sum_{i=0}^{n}(y_i - \hat{y}_i)^2}{\sum_{i=0}^{n}(y_i - \underline{y})^2} \ . \tag{9}$$

$(y_i$: actual value, $\hat{y}_i: predicted\ value$, $\underline{y}$: mean value)

# 3. Results and Analysis

## 3.1 Correlation Analysis of material parameters and indentation depth

Although an ANN can be designed using all eight material parameters ($A_1$, $A_2$, $B$, $n$, $C$, $\varepsilon$, $m$, and ISE), a large number of training datasets is required in eight dimensions. Thus, acquiring sufficient training data corresponding to eight parameters is quite challenging. To address this, we estimated the degree to which each material parameter affected the indentation depth before creating the training datasets. We termed this the absolute correlation coefficient.

For our correlation analysis, the eight material parameters were either multiplied or divided by a constant factor (i.e., 2 and 4) to generate 33 analytical datasets, as shown in **Table A.1**. Simulations were performed using these 33 datasets and the absolute correlation coefficients for each material parameter corresponding to different indentation depths were calculated, as shown in **Table 5.**

A larger absolute correlation coefficient for a given material parameter indicates that it has a greater effect on the indentation depth. In this study, we selected only those material parameters with an absolute correlation coefficient greater than 0.1. Based on this selection criterion, only the material parameters $B$, $n$, and $C$ of the J–C model qualified for all three impact cases. Consequently, the number of input variables in the ANN was reduced from eight to three. These three variables were multiplied and divided by factors of 2 and 4 to generate 125 training datasets for the ANN. **Table 6** shows the parameter settings for these training datasets.

## 3.2 Design and selection of optimal ANN models

The three material parameters *B*, *n*, and *C*, which had the most influence on the indentation depth, formed the input layer of our ANN, whereas the indentation depth formed the output layer. The details of the activation functions, learning rate, and epoch used in our model have already been discussed in **Section 2.2**. We named our models A, B, and C corresponding to impact cases 1, 2, and 3, respectively. To increase the prediction accuracy of the ANN having a single hidden layer, the total number of neurons was increased. These models were named A1, A2, A3, and A4 corresponding to 25, 50, 75, and 100 neurons, respectively.

**Figure. 5** shows the RMSE and $R^2$ as functions of the number of neurons and number of hidden layers for different ANN models. The smaller the value of RMSE and closer the value of $R^2$ to 1.0, the higher is the accuracy of the model. We observe that for a single hidden layer, the accuracy of the ANN was enhanced when the number of neurons was increased to 100 for all three impact cases.

In addition, we varied the number of hidden layers (i.e., 1, 2, and 3) of the ANN while keeping the total number of neurons fixed at 100. For example, model A4-2 represents impact case 1 and an ANN composed of 100 neurons and 2 hidden layers. Based on our evaluation, models A4-2, B4-3, and C4-2 were found to have the best prediction accuracy. This implies that increasing the number of hidden layers does not necessarily guarantee better results. Although the prediction accuracy improved as the number of hidden layers were increased in model B4, this was not the case for models A4 and C4. For example, the RMSE of model A4-2 increased from 0.4493 to 0.5106 for model A4-3 while $R^2$ remained constant at 0.999. Similarly, the RMSE of model C4-2 increased from 0.4222 to 0.4298 for model C4-3 while $R^2$ remained constant at 0.999.

## 3.3 Prediction of optimal material parameters

To find the optimal material parameters within the range of our training datasets, the maximum and minimum values of the three parameters *B*, *n*, and *C* were divided by a constant factor. The maximum and minimum values of these parameters are listed in **Table 7**. Subsequently, one million target datasets were generated from 100 variations of these three input parameters.

Based on our evaluation, models A4-2, B4-3, and C4-2 were tested on the one million target datasets to predict the indentation depth for each impact case. Next, the function *f* was calculated as the sum of the differences between the indentation depth predicted by the ANN models and the reference indentation depth for each impact case (**Standard, 1987; Standard, 2008; Kim, 2018**), such that

$$f = \sum_{j}^{N}\{(|16 - k_{1j}| + |25 - k_{2j}| + |19 - k_{2j}|)/3\} \ . \tag{10}$$

Finally, the optimal material parameters were selected using the inverse tracking method ($f^{-1}$).

**Figure. 6** shows the geometric distributions of the material parameters predicted by the ANN models and inverse tracking method with the following characteristics. In **Figure. 6**, 5,000 sets of material parameters, which induce low level of error less than 10 % between the indentation depths predicted by our algorithm and the corresponding reference indentation depths, are shown for each impact case. Meanwhile, the data highlighted in red correspond to 50 sets of material parameters that are common to all three impact cases with an averaged error below 11%.

In **Figure. 6**, we observe that the geometric distributions of the material parameters corresponding to impact cases 1 and 3 are almost indistinguishable; however, the geometric distribution corresponding to impact case 2 is distinctly different. Three separate intersection regions are also seen to form, which consist of 50 sets of material parameters.

### 3.4 Data validation

The ten best sets of material parameters that are applicable to all three impact cases were selected from among the 50 sets of material parameters discussed in **Section 3.3**. These ten sets are presented in **Table B.1.** Note that the average error between the indentation depths calculated using these ten sets of material parameters and the corresponding reference indentation depths is expected to be less than 7%; nevertheless, these material parameters need to be validated. Therefore, we conducted FEM simulations using these ten sets of material parameters to verify that the numerical indentation depths were in agreement with those predicted by the ANN models. The absolute error between the indentation depths yielded by the FEM simulations and those predicted by our ANN models for the ten best sets of material parameters was compared for all three impact cases, as shown in **Figure. 7**.

Following the data validation step, the optimal material parameters for modeling RP clay were selected, as listed in **Table 8**. The parameters $B$, $n$, and $C$ were determined using our ANN models, while the remaining parameters were adopted from previous studies. We compared the indentation depths that were numerically estimated using the optimal material parameters (predicted by our ANN models) with those that were estimated using the default material parameters, as shown in **Figure. 8**. **Figure. 9** displays the percentage relative error with respect to the reference indentation depth for the cases presented in **Figure. 8**. We observe that the use of optimal material parameters significantly reduces the mean relative error from 17.27% to 1.74%.

## 3.5 Effect of the number of training datasets

As mentioned earlier, a sufficient number of training datasets is required to ensure that the ANN has good prediction accuracy. However, because of its black-box nature, it is difficult to determine the exact number of training datasets required by an ANN. Therefore, we varied the number of training datasets to examine its effect on the prediction accuracy of the ANN. Our model was trained using a randomly reduced number of datasets (i.e., 75%, 50%, 25%, 10%, and 5%) from a set of 125 previously generated datasets.

When ANNs are trained using a randomly extracted number of training datasets, the reliability of a single result is difficult to guarantee. Therefore, we designed 10 models using the reduced training dataset for each impact case. For example, first, we trained an ANN model using 75% of the original 125 training datasets for impact case 1; next, this process was repeated 10 times to generate 10 ANN models for impact case 1.

The ANN models trained on a reduced number of datasets were used to predict the indentation depths corresponding to the input parameters $B$, $n$, and $C$. The reference indentation depth was compared with the indentation depth predicted using ten sets of material parameters having the lowest error that were selected using the inverse tracking method. The detailed process is illustrated in **Figure. 1**. FEM simulations were performed using the optimal material parameters, and the indentation depths predicted by the ANN models were compared with the simulation results. **Figure. 10** shows the percentage relative error between the FEM and ANN results as well as the average prediction accuracy of the ANN models for a reduced number of training datasets.

We observed that the average prediction accuracy of the ANN models declined as the number of training datasets was reduced. For example, the average accuracy was more than 90% when more than 25% of the training datasets were used; however, the average accuracy dropped sharply when less than 10% of the training datasets were used. In addition, the relative error remained similar when more than 25% of the training datasets were used. Thus, at least 31 training datasets are required to ensure a prediction accuracy of 90% or more for our ANN models.

Furthermore, a larger number of training datasets does not necessarily guarantee a higher prediction accuracy for different impact cases. For example, for impact case 2, the relative error of the model designed using 50% of the training datasets was 2% higher than that of the model designed using only 25% of the training dataset. This may be the result of a significant loss in specific training datasets during the random extraction from the original datasets. Nevertheless, the average accuracy of the ANN models, considering all three impact cases, increased as the number of training datasets were increased.

# 4. Discussion

In this study, we used ANNs as a tool to determine the optimal material parameters required to model RP clay for three different impact cases. In this section, we discuss some of the implications of our findings.

First, the parameters *B*, *n*, and *C* are directly related to the strain rate were most influential, among the eight material parameters of the RP clay model. This was confirmed by the correlation analysis performed for all three impact cases presented in **Section 3.1**. Considering the relatively low impact velocities of the cases and the correlation analysis result, it is mandatory to include the effect of strain rate to obtain more accurate numerical models.

Second, impact case 2 is the key to determine the optimal material parameters that are applicable to all three impact cases. The geometric distributions of the material parameters, which were selected within a certain error, were similar for the different impact cases except for impact case 2, represented as a blue plane in **Figure. 6**. This difference is responsible for the formation of complex intersection key regions between the material parameters of the three impact cases in **Figure. 6**. This distinctive geometric difference in **Figure. 6** seems to be caused by different strain rate in between the cases. The ratio between impact velocity and diameter of indenter was almost identical in impact cases 1 and 3. However, the ratio of impact case 2 was about 1.4 times higher. The higher ratio provides higher range of strain rate during deformation. This implies that the ANN models which dealt with strain rate effect need to have various impact cases having wide range of strain rate. The selection of comparable target dataset is crucial for determining the optimal material parameters. If there were no different strain rates in impact cases, the optimal material parameters or intersecting key region would be hardly determined.

Third, only 31 training datasets were sufficient to ensure a prediction accuracy of more than 90% for our ANN models. It is important to determine the optimal number of training datasets required to make accurate predictions using ANNs. We found that the prediction accuracy of our ANN models remained above 90% even when only 25% of the original number (125) of training datasets were used. However, when the data loss exceeded 75%, the prediction accuracy of the models decreased sharply. A high prediction accuracy of 90% even for a significantly reduced training dataset is probably because our ANNs consist of a large number of neurons and multiple hidden layers.

It is notable to observe that the level of accuracy was improved dramatically from 17.27% to 1.74% by incorporating FEM and ANN models. Our research work implies that the selection of appropriate cases covering different strain rate is important as discussed with **Figure 6.** On the other hand, the proposed simulation models considered only frictionless condition and indentation depth as comparable data. Any more informative series of comparable data, for example, geometric deformation contour near indentation, velocity profile of the indenter from experiment, frictional effect, and microfracture in the contacting layer may be helpful to understand and improve numerical model. The result shows that the hybrid method can also easily be used to simulate other impact cases with considering necessary type of data and key cases containing physical phenomenon, similarly to strain rate effect. The proposed method would relieve the burden of building numerical models, which is to get input parameters from tons of various experiments.

# 5. Conclusions

In this study, we determined the optimal values of the material parameters required to model RP clay, which is commonly used to evaluate bulletproof armors, using ANNs and the inverse tracking method. Among the eight material parameters (i.e. $A_1$, $A_2$, $B$, $n$, $C$, $\varepsilon$, $m$, and ISE) typically used to model RP clay, we selected only the strain-rate-dependent material parameters of the J–C model, namely $B$, $n$, and C. The optimal $B$, $n$, and $C$ values returned by our ANN models were 86.545 kPa, 0.171 kPa, and 0.479, respectively. We found that the mean relative error between the referred indentation depth and the indentation depth numerically estimated using the default parameter values was significantly reduced from 17.27% to 1.74% when the optimal parameter values were used. The geometric distribution with the parameter inputs emphasized that the range target dataset must be desirably selected to provide intersecting key region. Our proposed algorithm can be used to create an accurate clay model, which in turn can improve the effectiveness of bulletproof vests.

**Acknowledgements**

This work was supported by the Basic Science Research Program of the National Research Foundation (NRF) of Korea funded by the Ministry of Education (grant number: NRF-2021R1I1A3060343).

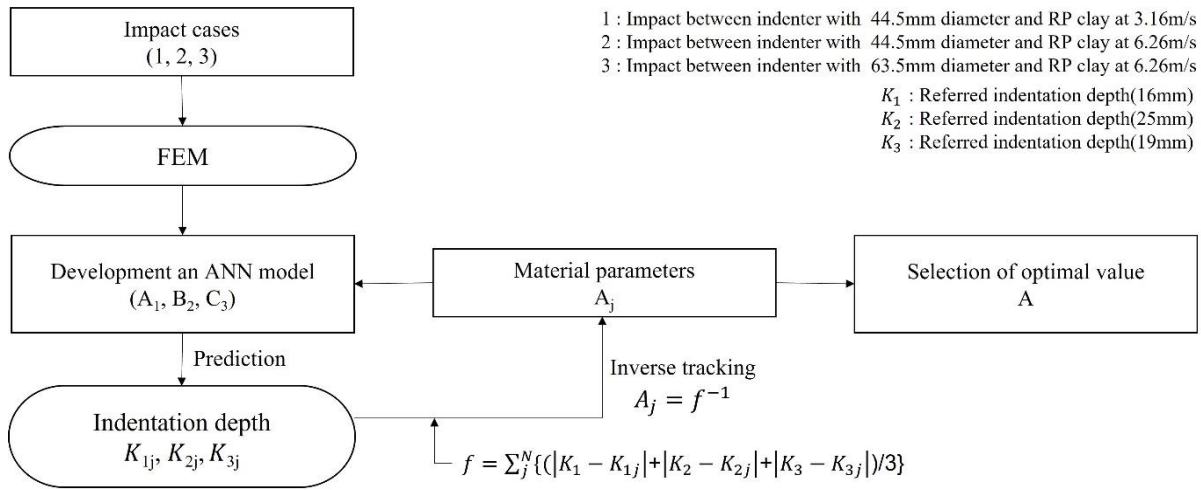

figure 1. Algorithm for determining the optimal material parameters of the RP clay model.

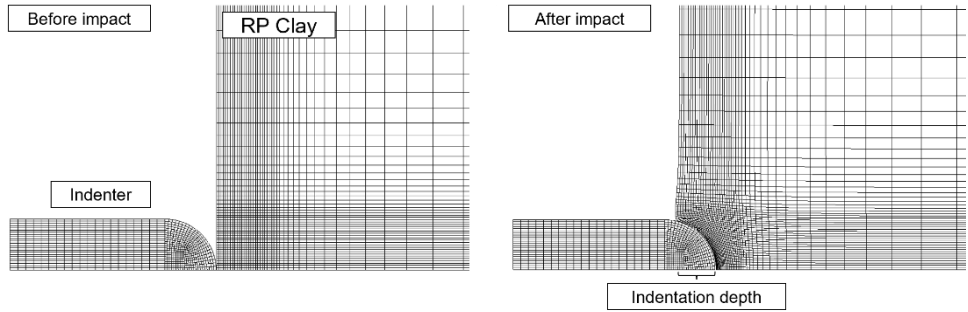

Figure 2. Schematic of the simulation setup for impact cases 1 and 2.

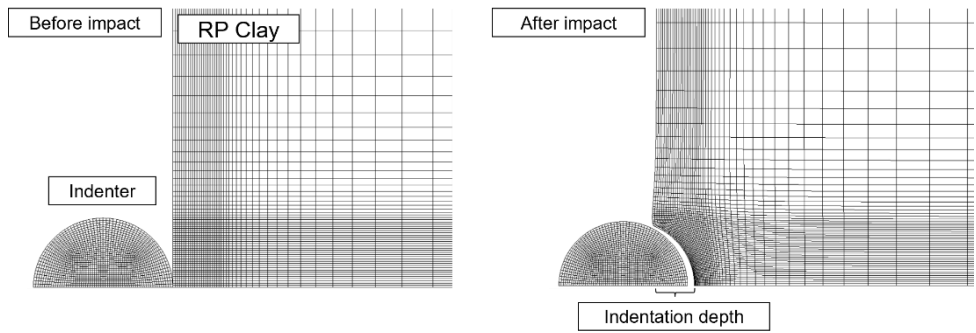

Figure 3. Schematic of the simulation setup for impact case 3.

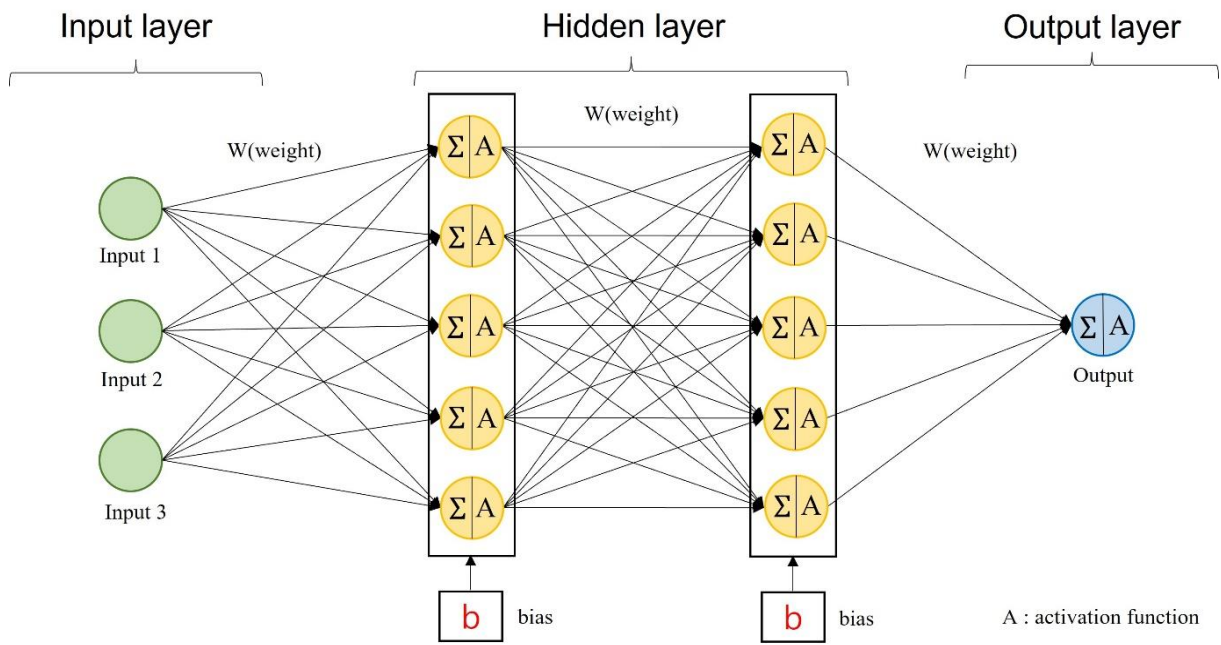

**Figure 4. Structure of a simplified ANN.**

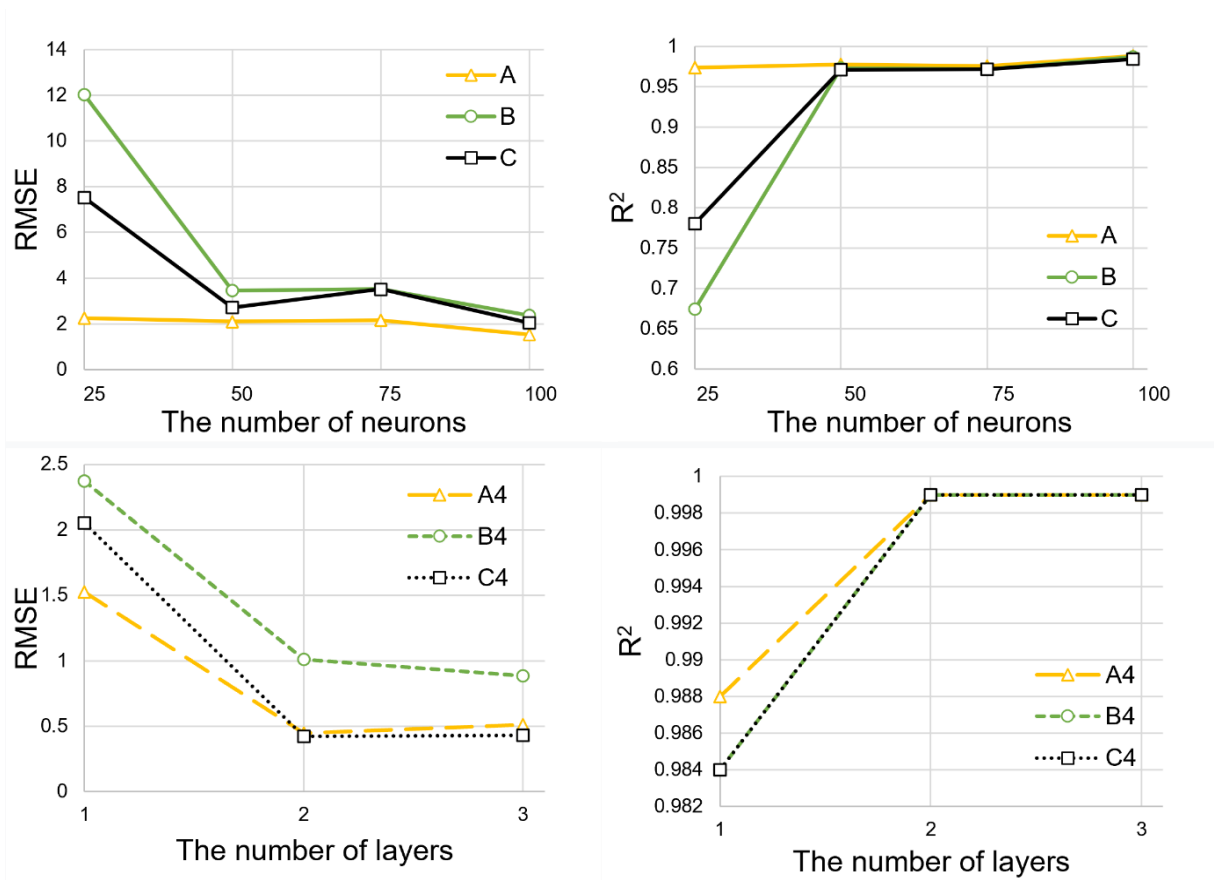

**Figure 5.** *Evaluation of the prediction accuracy of different ANN models.*

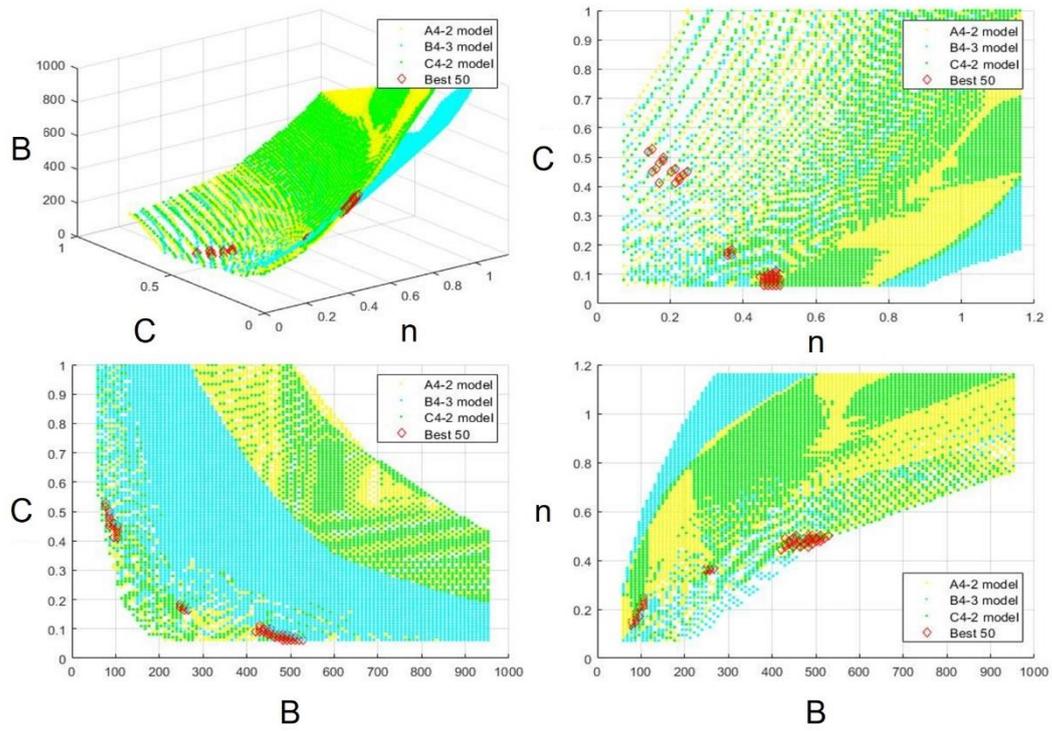

**Figure 6. Geometric distributions of the predicted material parameters.**

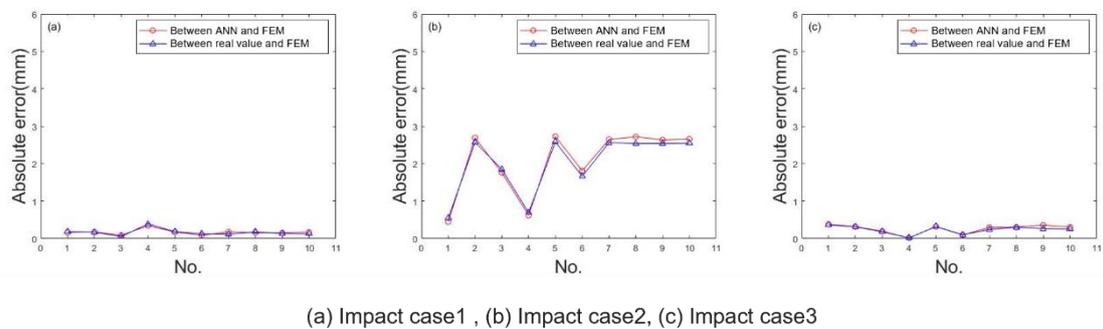

(a) Impact case1 , (b) Impact case2, (c) Impact case3

**Figure 7. Comparison of the FEM and ANN results for the ten best sets of material parameters predicted by the ANN models.**

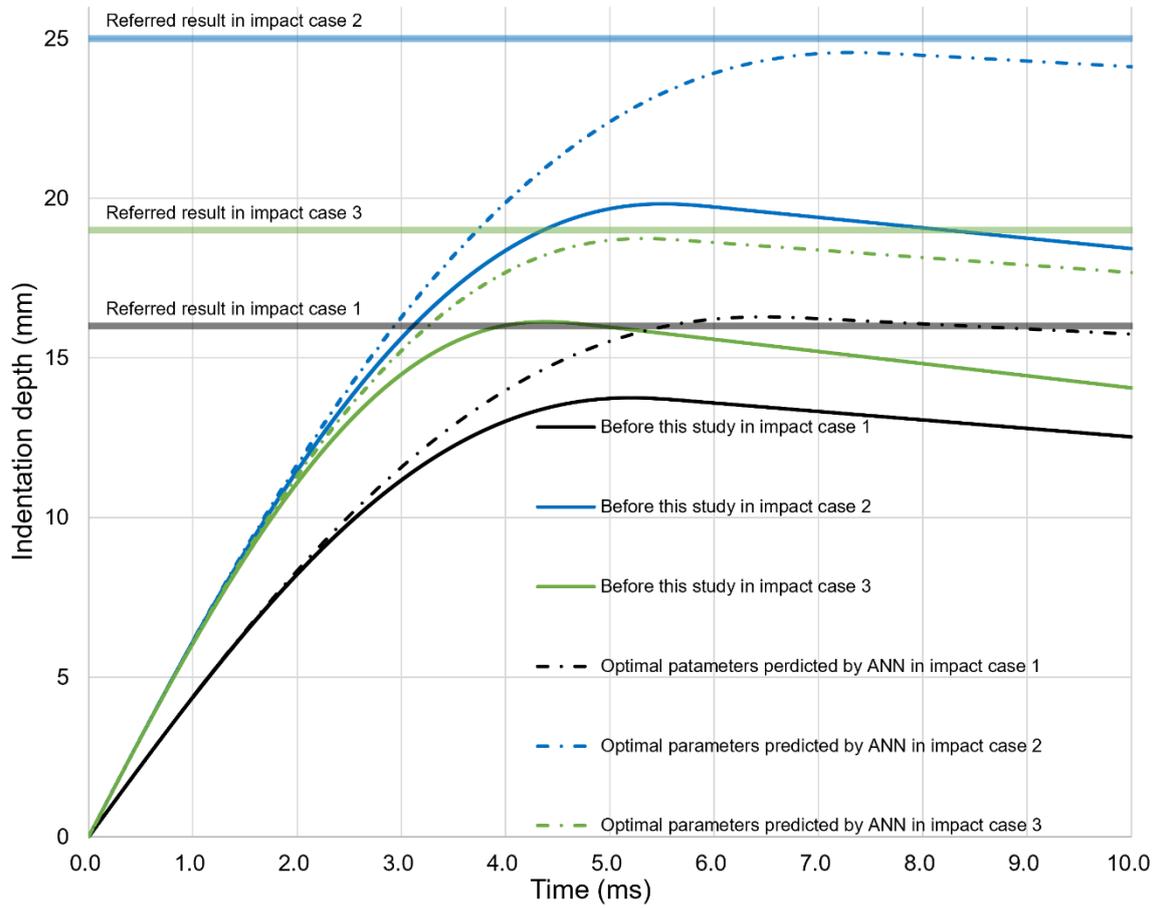

**Figure 8. Comparison of the indentation depths before and after using the optimal material parameters for different impact cases.**

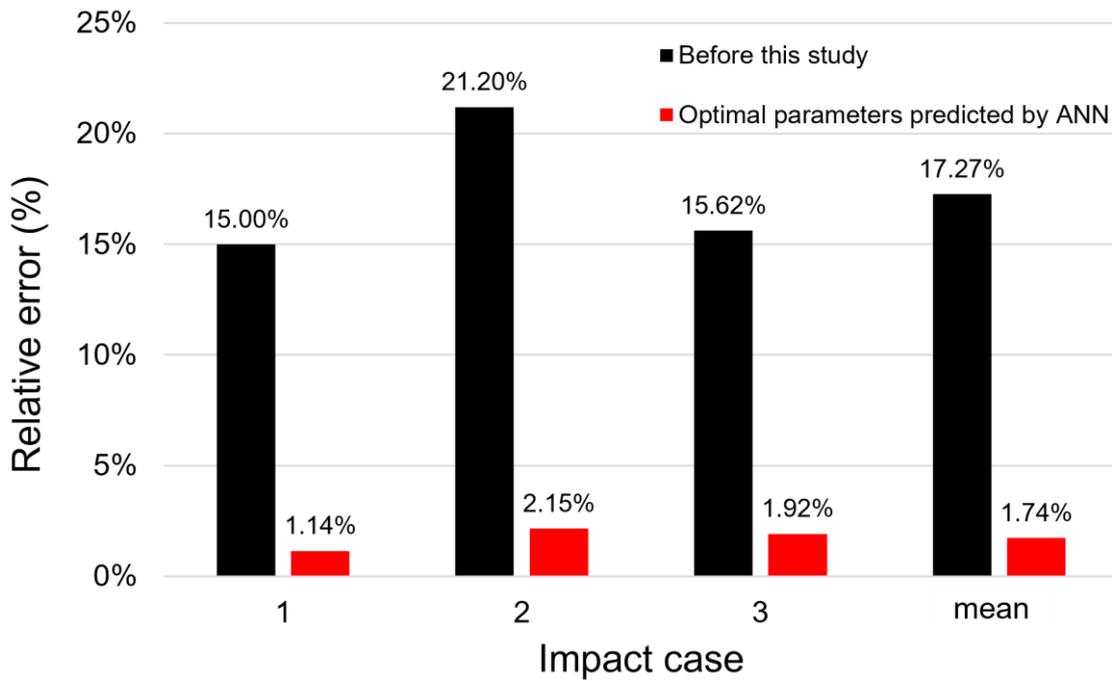

**Figure 9. Comparison of the percentage relative errors with respect to the reference indentation depth for the cases.**

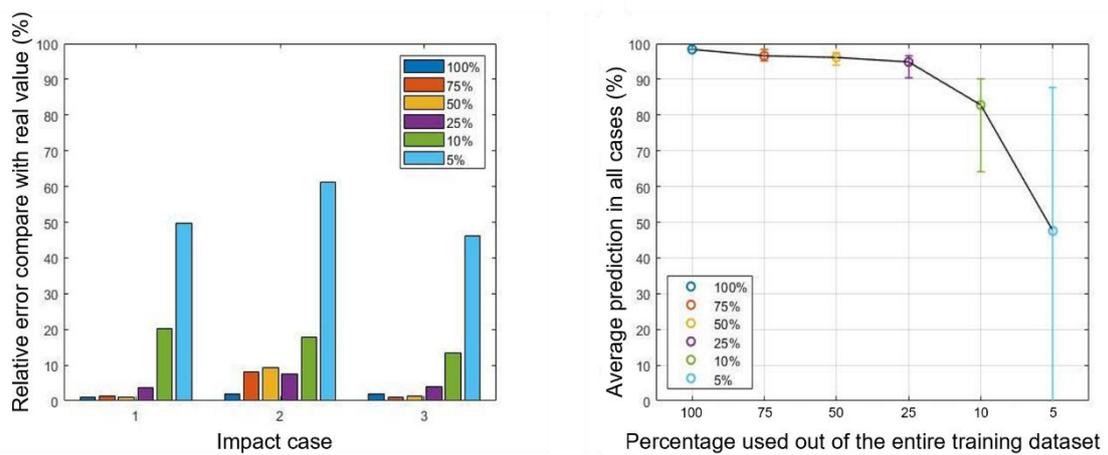

**Figure 10. Percentage relative error for each impact case (left) and average prediction accuracy (right) for a reduced number of training datasets.**

| Case | Diameter (mm) | Mass (g) | Impact velocity (m/s) | Indentation depth (mm) |
|---|---|---|---|---|
| 1 | 44.5 ±0.5 | 1,000±10 | 4.47 | 16.0 |
| 2 | 44.5 ±0.5 | 1,000±10 | 6.26 | 25.0 ± 3.0 |
| 3 | 63.5 ±0.05 | 1,043±5 | 6.26 | 19.0 ± 2.0 |

*Table 1. Parameter settings and experimental results of different impact cases (Standard, 1987; Standard, 2008; Kim, 2018).*

| Clay | | | Indenter | |
|---|---|---|---|---|
| Size (mm*mm) | Mesh size (mm) | Dimension ratio | Mesh size (mm) | |
| 400*200 | 1.0 | 10.0 | 1.0 | |

*Table 2. Mesh size and boundary conditions of both clay and indenter (Kim et al., 2022).*

| Parameter | $A_1$ (GPa) | $A_2$ (GPa) | $A_3$ (GPa) | $B_0$ | $B_1$ | |
|---|---|---|---|---|---|---|
| Value | 2.804 | 40.7 | -36.0 | 1.7 | 1.7 | |
| Parameter | A (kPa) | B (kPa) | n | C | m | $\varepsilon$ ($s^{-1}$) |
| Value | 0.01 | 238.0 | 0.29 | 0.25 | 0.502 | 0.118 |

*Table 3. Parameters of the polynomial EOS and J–C model for RP clay (Mates et al., 2014; Gad and Gao, 2020).*

| EOS | | J–C model | | | | | | Erosion |
|---|---|---|---|---|---|---|---|---|
| $A_1$ (GPa) | $A_2$ (GPa) | B (kPa) | n | C | m | $\varepsilon$ ($s^{-1}$) | | ISE |
| 2.804 | 40.7 | 238.0 | 0.29 | 0.25 | 0.502 | 0.118 | | 5.0 |

*Table 4. Initially selected parameters for modeling RP clay.*

| Material parameters | Absolute correlation coefficient | | |
|---|---|---|---|
| | Impact cases | | |
| | 1 | 2 | 3 |
| A1 | 0.029 | 0.031 | 0.030 |
| A2 | 0.029 | 0.031 | 0.030 |
| B | 0.522 | 0.500 | 0.470 |
| n | 0.600 | 0.490 | 0.636 |
| C | 0.289 | 0.309 | 0.297 |
| $\varepsilon$ | 0.025 | 0.026 | 0.024 |
| m | 0.070 | 0.078 | 0.074 |
| ISE | 0.029 | 0.034 | 0.030 |

*Table 5. Absolute correlation coefficients for different indentation depths.*

| Variable material parameters | | | Default material parameters | | | | |
|---|---|---|---|---|---|---|---|
| B (kPa) | n | C | A1 (GPa) | A2 (GPa) | m | $\varepsilon$ | ISE |
| 59.5 | 0.0725 | 0.0625 | | | | | |
| 119.0 | 0.1450 | 0.1250 | | | | | |
| 238.0 | 0.2900 | 0.2500 | 2.804 | 40.7 | 0.502 | 0.118 | 5.0 |
| 476.0 | 0.5800 | 0.5000 | | | | | |
| 952.0 | 1.1600 | 1.0000 | | | | | |

*Table 6. Parameter settings of the training datasets.*

| | Material parameters | | |
|---|---|---|---|
| | B (kPa) | n | C |
| Minimum value | 59.5 | 0.0725 | 0.0625 |
| Maximum value | 952.0 | 1.1600 | 1.0000 |

*Table 7. Range of variable parameters for the one million target datasets.*

| J–C model | | | | | | EOS model | | | | |
|---|---|---|---|---|---|---|---|---|---|---|
| A (kPa) | B (kPa) | n | C | m | $\varepsilon$ ($s^{-1}$) | A1 (GPa) | A2 (GPa) | A3 (GPa) | B0 | B1 |
| 0.01 | 86.545 | 0.171 | 0.479 | 0.502 | 0.118 | 2.804 | 40.70 | -36.0 | 1.70 | 1.70 |

*Table 8. Optimal material parameters of the EOS and J–C model for RP clay.*

*Appendix A.*

| No. | Material parameters | | | | | | | Erosion | Indentation depth (mm) | | |
|---|---|---|---|---|---|---|---|---|---|---|---|
| | EOS model | | J–C model | | | | | | Impact cases | | |
| | A1 (GPa) | A2 (GPa) | B (kPa) | n | C | $\varepsilon$ $(s^{-1})$ | m | ISE | 1 | 2 | 3 |
| 1 | 0.701 | 40.700 | 238.0 | 0.2900 | 0.2500 | 0.1180 | 0.5020 | 5.00 | 13.550 | 19.680 | 16.020 |
| 2 | 1.402 | 40.700 | 238.0 | 0.2900 | 0.2500 | 0.1180 | 0.5020 | 5.00 | 13.560 | 19.720 | 16.030 |
| 3 | 5.608 | 40.700 | 238.0 | 0.2900 | 0.2500 | 0.1180 | 0.5020 | 5.00 | 13.565 | 19.750 | 16.035 |
| 4 | 11.216 | 40.700 | 238.0 | 0.2900 | 0.2500 | 0.1180 | 0.5020 | 5.00 | 13.557 | 19.700 | 16.027 |
| 5 | 2.804 | 10.175 | 238.0 | 0.2900 | 0.2500 | 0.1180 | 0.5020 | 5.00 | 13.563 | 19.740 | 16.033 |
| 6 | 2.804 | 20.350 | 238.0 | 0.2900 | 0.2500 | 0.1180 | 0.5020 | 5.00 | 13.563 | 19.740 | 16.033 |
| 7 | 2.804 | 81.400 | 238.0 | 0.2900 | 0.2500 | 0.1180 | 0.5020 | 5.00 | 13.563 | 19.740 | 16.033 |
| 8 | 2.804 | 162.800 | 238.0 | 0.2900 | 0.2500 | 0.1180 | 0.5020 | 5.00 | 13.563 | 19.740 | 16.033 |
| 9 | 2.804 | 40.700 | 59.5 | 0.2900 | 0.2500 | 0.1180 | 0.5020 | 5.00 | 31.525 | 49.840 | 35.384 |
| 10 | 2.804 | 40.700 | 119.0 | 0.2900 | 0.2500 | 0.1180 | 0.5020 | 5.00 | 20.110 | 30.500 | 23.262 |
| 11 | 2.804 | 40.700 | 476.0 | 0.2900 | 0.2500 | 0.1180 | 0.5020 | 5.00 | 9.500 | 13.480 | 11.415 |
| 12 | 2.804 | 40.700 | 952.0 | 0.2900 | 0.2500 | 0.1180 | 0.5020 | 5.00 | 5.030 | 9.620 | 8.363 |
| 13 | 2.804 | 40.700 | 238.0 | 0.0725 | 0.2500 | 0.1180 | 0.5020 | 5.00 | 10.000 | 14.950 | 11.777 |
| 14 | 2.804 | 40.700 | 238.0 | 0.1450 | 0.2500 | 0.1180 | 0.5020 | 5.00 | 11.074 | 16.410 | 13.064 |
| 15 | 2.804 | 40.700 | 238.0 | 0.5800 | 0.2500 | 0.1180 | 0.5020 | 5.00 | 19.330 | 27.060 | 22.790 |
| 16 | 2.804 | 40.700 | 238.0 | 1.1600 | 0.2500 | 0.1180 | 0.5020 | 5.00 | 30.100 | 39.450 | 35.330 |
| 17 | 2.804 | 40.700 | 238.0 | 0.2900 | 0.0625 | 0.1180 | 0.5020 | 5.00 | 18.800 | 28.150 | 21.850 |
| 18 | 2.804 | 40.700 | 238.0 | 0.2900 | 0.1250 | 0.1180 | 0.5020 | 5.00 | 16.420 | 24.450 | 19.210 |
| 19 | 2.804 | 40.700 | 238.0 | 0.2900 | 0.5000 | 0.1180 | 0.5020 | 5.00 | 10.680 | 15.200 | 12.775 |
| 20 | 2.804 | 40.700 | 238.0 | 0.2900 | 1.0000 | 0.1180 | 0.5020 | 5.00 | 8.155 | 11.400 | 9.860 |
| 21 | 2.804 | 40.700 | 238.0 | 0.2900 | 0.2500 | 0.0295 | 0.5020 | 5.00 | 12.620 | 18.300 | 14.98 |
| 22 | 2.804 | 40.700 | 238.0 | 0.2900 | 0.2500 | 0.0590 | 0.5020 | 5.00 | 13.065 | 18.980 | 15.478 |
| 23 | 2.804 | 40.700 | 238.0 | 0.2900 | 0.2500 | 0.2360 | 0.5020 | 5.00 | 14.125 | 20.600 | 16.660 |
| 24 | 2.804 | 40.700 | 238.0 | 0.2900 | 0.2500 | 0.4720 | 0.5020 | 5.00 | 14.760 | 21.600 | 17.368 |
| 25 | 2.804 | 40.700 | 238.0 | 0.2900 | 0.2500 | 0.1180 | 0.5020 | 1.25 | 13.563 | 20.380 | 16.033 |
| 26 | 2.804 | 40.700 | 238.0 | 0.2900 | 0.2500 | 0.1180 | 0.5020 | 2.500 | 13.563 | 19.740 | 16.033 |
| 27 | 2.804 | 40.700 | 238.0 | 0.2900 | 0.2500 | 0.1180 | 0.5020 | 10.00 | 13.563 | 19.740 | 16.033 |
| 28 | 2.804 | 40.700 | 238.0 | 0.2900 | 0.2500 | 0.1180 | 0.5020 | 20.00 | 13.563 | 19.740 | 16.033 |
| 29 | 2.804 | 40.700 | 238.0 | 0.2900 | 0.2500 | 0.1180 | 0.1255 | 5.00 | 17.468 | 26.08 | 20.220 |
| 30 | 2.804 | 40.700 | 238.0 | 0.2900 | 0.2500 | 0.1180 | 0.2510 | 5.00 | 14.700 | 20.580 | 17.200 |
| 31 | 2.804 | 40.700 | 238.0 | 0.2900 | 0.2500 | 0.1180 | 1.0040 | 5.00 | 13.476 | 19.400 | 15.850 |
| 32 | 2.804 | 40.700 | 238.0 | 0.2900 | 0.2500 | 0.1180 | 2.0080 | 5.00 | 13.470 | 19.400 | 15.840 |
| 33 | 2.804 | 40.700 | 238.0 | 0.2900 | 0.2500 | 0.1180 | 0.5020 | 5.00 | 13.600 | 19.700 | 16.033 |

*Table A.1. Details of the analysis datasets.*

*Appendix B.*

| No. | Material parameters | | | Indentation depth (mm) | | | | | |
| --- | --- | --- | --- | --- | --- | --- | --- | --- | --- |
| | | | | ANN Impact cases | | | FEM Impact cases | | |
| | B (kPa) | n | C | 1 | 2 | 3 | 1 | 2 | 3 |
| 1 | 86.545 | 0.171 | 0.4792 | 16.020 | 24.909 | 19.014 | 16.183 | 24.463 | 18.636 |
| 2 | 456.167 | 0.468 | 0.0814 | 16.005 | 25.120 | 19.007 | 15.830 | 22.430 | 18.690 |
| 3 | 257.833 | 0.369 | 0.1761 | 15.965 | 24.910 | 19.014 | 16.050 | 23.150 | 18.820 |
| 4 | 104.576 | 0.237 | 0.4419 | 16.047 | 24.930 | 18.971 | 16.390 | 24.310 | 18.985 |
| 5 | 474.197 | 0.468 | 0.0720 | 15.993 | 25.147 | 18.996 | 15.820 | 22.420 | 18.680 |
| 6 | 248.818 | 0.358 | 0.1761 | 16.037 | 25.149 | 18.993 | 16.125 | 23.340 | 18.900 |
| 7 | 465.182 | 0.479 | 0.0814 | 16.058 | 25.083 | 19.055 | 15.885 | 22.440 | 18.760 |
| 8 | 492.227 | 0.468 | 0.0625 | 15.981 | 25.178 | 19.003 | 15.825 | 22.460 | 18.700 |
| 9 | 438.136 | 0.468 | 0.0909 | 16.018 | 25.093 | 19.091 | 15.868 | 22.460 | 18.735 |
| 10 | 483.212 | 0.479 | 0.0720 | 16.046 | 25.112 | 19.053 | 15.876 | 22.455 | 18.750 |

*Table B.1. Comparison of the FEM and ANN indentation depths for the ten best sets of material parameters predicted by the ANN models.*